\documentclass[11pt, a4paper]{amsart}

\usepackage{amsthm}
\usepackage{url}

\newcommand{\rarr}{\rightarrow}
\newcommand{\sm}{\setminus}

\newcommand{\R}{\mathbb{R}}
\newcommand{\N}{\mathbb{N}}

\newcommand{\mP}{\mathcal{P}}

\newtheorem{thm}{Theorem}
\newtheorem{lemma}[thm]{Lemma}

\begin{document}

\title{Linear Bound for Majority Colourings of Digraphs}
\author{Fiachra Knox}
\address[Fiachra Knox]{Department of Mathematics\\
Simon Fraser University\\
Burnaby, British Columbia, Canada}
\email{fknox@sfu.ca}
\thanks{Research of the first author is supported by a PIMS postdoctoral fellowship.}

\author{Robert \v{S}\'{a}mal} 
\address[Robert \v{S}\'{a}mal]{Computer Science Institute\\ Charles University\\ Prague, Czech Republic} 
\email{samal@iuuk.mff.cuni.cz}
\thanks{The second author is partially supported by grant GA \v{C}R 16-19910S and partially supported by grant LL1201 ERC CZ of the Czech Ministry of Education, Youth and Sports.}

\date{\today}

\maketitle

\begin{abstract} Given $\eta \in [0, 1]$, a colouring $C$ of $V(G)$ is an \emph{$\eta$-majority colouring} if at most $\eta d^+(v)$ out-neighbours of $v$ have colour $C(v)$, for any $v \in V(G)$.
We show that every digraph $G$ equipped with an assignment of lists $L$, each of size at least $k$, has a $2/k$-majority $L$-colouring.
For even $k$ this is best possible, while for odd $k$ the constant $2/k$ cannot be replaced by any number less than $2/(k+1)$.
This generalizes a result of Anholcer, Bosek and Grytczuk~\cite{ABGpp}, who proved the cases $k=3$ and $k=4$ and gave a weaker result for general $k$.
\end{abstract}

\section{Introduction}

Given a digraph $G$, we write $V(G)$ and $E(G)$ for the vertex and edge set of a digraph $G$, respectively.
For $v \in V(G)$, we denote by $d^+(v)$ the out-degree of $v$.
Given $\eta \in [0, 1]$, a (not necessarily proper) colouring $C$ of $V(G)$ is an \emph{$\eta$-majority colouring} if at most $\eta d^+(v)$ out-neighbours of $v$ have colour $C(v)$, for any $v \in V(G)$.
A $1/2$-majority colouring is referred to simply as a \emph{majority colouring}.
This concept was introduced by van der Zypen~\cite{vdZmathOv}, who asked whether every digraph has a majority colouring with a bounded number of colours.
This question was answered by Kreutzer, Oum, Seymour, van der Zypen and Wood~\cite{ZypenWoodetal}, who showed that $4$ colours always suffice.

We consider the list-colouring version of this problem. 
For a set $S$, we denote by $\mP(S)$ the power set of $S$.
Given a digraph $G$ and an assignment $L: V(G) \rarr \mP(\N)$ of lists to vertices of $G$, an \emph{$L$-colouring} $C: V(G) \rarr \N$ of $G$ is a colouring of $V(G)$ such that $C(v) \in L(v)$ for every $v \in V(G)$.
If $G$ has an $\eta$-majority $L$-colouring for any such assignment $L$ whose lists are all of size at least $k$, we say that $G$ is \emph{$\eta$-majority $k$-choosable}.
Anholcer, Bosek and Grytczuk~\cite{ABGpp} showed that every graph $G$ is $1/n$-majority $n^2$-choosable for every $n \geq 2$.
Theorem~\ref{thm:Main} improves on this result.

\begin{thm} \label{thm:Main} For any integer $k \geq 2$, every digraph $G$ is $2/k$-majority $k$-choosable.
\end{thm}

Theorem~\ref{thm:Main} was proved independently by Gir\~{a}o, Kittipassorn and Popielarz~\cite{GKPpp}.
The case $k=2$ is trivial.
Previously, Anholcer, Bosek and Grytczuk~\cite{ABGpp} showed that Theorem~\ref{thm:Main} holds in the cases $k=3$ and $k=4$ and conjectured that $2/k$ can be replaced by $1/2$ when $k = 3$.
Theorem~\ref{thm:Main} is best possible when $k$ is even, as shown by the example of a $k/2$-regular tournament on $k+1$ vertices (that is, all vertices 
have both in-degree and out-degree equal to $k/2$). If we make all lists equal, 
then some vertex must have an out-neighbour of the same colour, and this out-neighbour represents $2/k$ of its out-neighbourhood.
When $k$ is odd, a similar example shows that we cannot replace $2/k$ by any number less than $2/(k+1)$.

\section{Proof of Theorem~\ref{thm:Main}}

We denote by $vw$ an edge from a vertex $v$ of a digraph to another vertex~$w$.
The proof of Theorem~\ref{thm:Main} relies on Lemma~\ref{lem:Key}, which follows.

\begin{lemma} \label{lem:Key} Let $k \geq 2$ be an integer and let $G$ be a digraph on a vertex set $V = S \cup T$, such that $G[S]$ is strongly connected, $G[T]$ is
  edgeless and there are no edges from $T$ to $S$.
Let $C_T$ be any colouring of $T$ and let $L:S\rarr \mP(\N)$ be an assignment of lists, each of size at least $k$, to vertices in $S$.
Then there exists an extension $C$ of $C_T$ to $V$ with $C(v) \in L(v)$ for each $v \in S$, such that no vertex $v \in S$ has more than $2d^+(v)/k$ out-neighbours with the same colour as $v$.
\end{lemma}

\proof For any colouring $C$ of $V$, we define the function $f_C: S \rarr \R$ by
$$f_C(v) = \frac{|\{w \in N^+(v) \mid C(w) = C(v)\}|}{d^+(v)}$$
for each vertex $v \in S$; 
i.e., $f_C(v)$ is the proportion of out-neighbours of $v$ which have the same colour as $v$ under $C$.
Given $v \in S$, we write $d^+_S(v) = |N^+(v) \cap S|$.

Let $A$ be the non-negative real $S \times S$ matrix with entries $A_{vw} = 1/d^+_S(v)$ if $vw$ is an edge of $G$ and $A_{vw} = 0$ otherwise.
We have $A{\bf j} = {\bf j}$ (where ${\bf j}$ is the vector of all $1$'s).
On the other hand if $A{\bf y} = c{\bf y}$ for any vector ${\bf y}$, then choose $v \in S$ such that $|y_v|$ is maximal;
now $|cy_v| = |\sum_{w \in S} A_{vw} y_w| \leq \sum_{w \in S} A_{vw} |y_v| = |y_v|$ and so $|c| \leq 1$.
Thus, the spectral radius of $A$ is $1$.

By applying the Perron--Frobenius Theorem (see, e.g., \cite[Theorem 8.8.1]{godsil2013algebraic}) to~$A^\top$, 
noting that $G[S]$ is strongly connected, we obtain an eigenvector ${\bf x}$ of~$A^\top$ with positive entries and eigenvalue $1$.
We remark that by normalizing~${\bf x}$ we could obtain a stationary distribution of the uniform random walk on $G[S]$.

Consider an extension $C$ of $C_T$ with $C(v) \in L(v)$ for each $v \in S$ such that $\sum_{v \in S} x_v f_C(v)$ is minimized.
We claim that $C$ satisfies the requirements of the lemma.
For brevity we write $f$ for $f_C$.
It suffices to show that $f(v) \leq 2/k$ for every $v \in S$.
Observe that 
\begin{equation} \label{eq:minimized sum}
\sum_{v \in S} x_v f(v) = \sum_{\stackrel{vw \in E(G)}{C(v) = C(w)}} \frac{x_v}{d^+(v)}.
\end{equation}

Fix a vertex $v \in S$.
We define $g: L(v) \rarr \R$ by
$$g(i) = \sum_{\stackrel{w \in N^+(v)}{C(w) = i}} \frac{x_v}{d^+(v)} + \sum_{\stackrel{u \in N^-(v)}{C(u) = i}} \frac{x_u}{d^+(u)}$$
for $i \in L(v)$.
Observe that if $v$ were recoloured with colour $i$, then (\ref{eq:minimized sum}) would change by $g(i) - g(C(V))$.
By the minimality of $C$ and the definition of $g$ we have that $g(i) \geq g(C(v)) \geq x_v f(v)$.
Since $A^{\top} {\bf x} = {\bf x}$,
$$x_v = \sum_{u \in N^-(v)} \frac{x_u}{d^+_S(u)} \geq \sum_{u \in N^-(v)} \frac{x_u}{d^+(u)}$$
and hence
$$2x_v \geq \sum_{i \in L(v)} g(i) \geq kx_v f(v).$$ 
Since $x_v > 0$, we have $f(v) \leq 2/k$.
It follows immediately that $C$ satisfies the requirements of the lemma.
\endproof

\proof[Proof of Theorem~\ref{thm:Main}] We partition $V(G)$ into strongly connected components $S_1, S_2, \ldots, S_r$, where there are no edges from $S_i$ to $S_j$ for any $i < j$.
Let $L:~V(G)\rarr \mP(\N)$ be an assignment of lists, each of size at least $k$, to vertices of $G$.
We write $A_i$ for $\bigcup_{j \in [i]} S_j$ (taking $A_0 = \emptyset$); 
let $C_0$ be the unique colouring of $A_0$.
For each $i = 0, 1, 2, \ldots, r-1$ in turn, we apply Lemma~\ref{lem:Key} with $S = S_{i+1}$, $T = A_i$ and $G = G[A_{i+1}] \sm G[A_i]$ to obtain an extension of $C_i$ to an $L$-colouring $C_{i+1}$ of $A_{i+1}$ 
such that no $v \in S_{i+1}$ (and hence no $v \in A_{i+1}$) has more than $2d^+(v)/k$ out-neighbours of the same colour.
At the end of this process we obtain $C_r$, which is the desired $2/k$-majority $L$-colouring of~$V(G)$.
\endproof

\section*{Acknowledgements} 

The authors would like to thank David Wood for presenting this problem at the workshop ``New Trends in Graph Colouring'' in Banff, October 2016, and for many 
helpful comments. 
We are also grateful to the organizers of this workshop and to our colleagues at the workshop for helpful discussions.

\bibliography{Majority}{}
\bibliographystyle{plain}

\end{document}